\documentclass[runningheads]{llncs}
\usepackage[T1]{fontenc}
\usepackage{graphicx}
\usepackage{booktabs}
\usepackage{color}

\usepackage{orcidlink}

\usepackage{amsmath}
\usepackage{subcaption}
\usepackage{float}
\spnewtheorem{assumption}[theorem]{Assumption}{\bfseries}{\itshape}
\spnewtheorem{application}{Application}{\bfseries}{\rmfamily}

\newcommand{\calH}{\mathcal{H}}

\newcommand{\calU}{\mathcal{U}}
\newcommand{\calX}{\mathcal{X}}
\newcommand{\calZ}{\mathcal{Z}}

\newcommand{\sbt}{\mathrm{subject\; to }}

\newcommand{\argmax}{\mathop{\rm argmax}}
\usepackage[misc]{ifsym}
\begin{document}
\title{Self-Supervised Penalty-Based Learning for Robust Constrained Optimization}

\titlerunning{Self-Supervised Learning for Robust Optimization}
\author{Wyame Benslimane\textsuperscript{\Letter}\orcidlink{0009-0006-7742-2932} \and Paul Grigas\orcidlink{0000-0002-5617-1058}}

\authorrunning{W. Benslimane and P. Grigas}
\institute{Department of Industrial Engineering and Operations Research \\ University of California-Berkeley, Berkeley, CA, USA  \\ \email{\{wyame.benslimane,pgrigas\}@berkeley.edu}}
\maketitle

\begin{abstract}
We propose a new methodology for parameterized constrained robust optimization, an important class of optimization problems under uncertainty, based on learning with a self-supervised penalty-based loss function. Whereas supervised learning requires pre-solved instances for training, our approach leverages a custom loss function derived from the exact penalty method in optimization to learn an approximation, typically defined by a neural network model, of the parameterized optimal solution mapping. Additionally, we adapt our approach to robust constrained combinatorial optimization problems by incorporating a surrogate linear cost over mixed integer domains, and a smooth approximations thereof, into the final layer of the network architecture.
We perform computational experiments to test our approach on three different applications: multidimensional knapsack with continuous variables, combinatorial multidimensional knapsack with discrete variables, and an inventory management problem. Our results demonstrate that our self-supervised approach is able to effectively learn neural network approximations whose inference time is significantly smaller than the computation time of traditional solvers for this class of robust optimization problems. Furthermore, our results demonstrate that by varying the penalty parameter we are able to effectively balance the trade-off between sub-optimality and robust feasibility of the obtained solutions.
\end{abstract}

\section{Introduction}
\vspace{-0.1cm}
Decision-making under uncertainty is an essential aspect of various real-world applications, ranging from facility location planning, portfolio optimization, and inventory management to healthcare and logistics. Robust optimization (RO) \cite{ben2009robust} offers a powerful and computationally viable framework for addressing decision-making under uncertainty. Most approaches used for robust optimization rely on constructing a tractable reformulation based on a carefully designed uncertainty set \cite{ben1999robust,ben2015deriving}, which is later fed into a standard optimization solver. Thus, solving robust optimization problems often involves specialized optimization algorithms like interior-point methods or other iterative algorithms that rely on first or second order information \cite{ben2009robust,bertsimas2022robustadaptiveoptim}. Although these methods have been well studied in theory, with strong convergence properties and theoretical guarantees, and applied in a variety of contexts, they face several limitations in terms of ease of computational efficiency and scalability to very large problems. These computational challenges are amplified even further when addressing demanding problems involving nonlinear functions, combinatorial constraints, and large-scale instances. The need for scalable robust optimization solutions, computed very efficiently with low latency, also stems from the need to solve real-time decision-making problems under time constraints in various applications. For instance, power systems require real-time adjustments to accommodate the fluctuating nature of renewable energy sources \cite{zamzam2019learningoptimalsolutionsextremely}, while ride-sharing platforms need to efficiently dispatch drivers to meet dynamic customer demand and traffic conditions \cite{kim2005optimal}.

On the other hand, the use of machine learning for solving optimization problem instances has shown great potential to address some of the above challenges. Several models have been developed to combine the strength of data-driven approaches to solve decision-making problems efficiently \cite{kotary2021endtoendconstrainedoptimizationlearning,amos2023tutorialamortizedoptimization}. While initial work in this area focused on improving the performance of unconstrained optimization \cite{li2016learning,andrychowicz2016learning,chen2022learning}, recent models have included methods for handling more general optimization problems \cite{kotary2021endtoendconstrainedoptimizationlearning}. For example, one line of work is to integrate constraints as a layer of a neural network using implicit differentiation  \cite{amos2021optnetdifferentiableoptimizationlayer,agrawal2019differentiableconvexoptimizationlayers,ferber2019mipaalmixedintegerprogram,vlastelica2020differentiationblackboxcombinatorialsolvers}. Another line of work consists of using a surrogate model relying on custom loss functions to train models that approximate continuous constrained optimization solutions \cite{park2022selfsupervisedprimalduallearningconstrained,donti2021dc3learningmethodoptimization,kotary2023predictthenoptimizeproxylearningjoint}. Still, the constrained problem formulations considered so far in this line of work have placed little focus on uncertainty, particularly as modeled with robust optimization.

While data-driven approaches have been used for robust optimization  \cite{bertsimas2014datadrivenrobustoptimization,wang2024learningdecisionfocuseduncertaintysets}, most approaches focus on learning better uncertainty set representations. Our work adopts a different angle by directly learning the solution to the RO problem using the learning to optimize framework. By approximating the solution to a class of parametric robust optimization problems, we do not have to rely on classical iterative solvers that can be computationally expensive, especially for time-sensitive and large-scale problems. Furthermore, we rely on a self-supervised \cite{park2022selfsupervisedprimalduallearningconstrained} penalty approach where we design a custom loss function based on the exact penalty method \cite{di1989exact} to construct a neural network-based solver for the RO problem. This approach allows us to potentially achieve faster solution times and improved performance compared to traditional methods. Furthermore, to handle a wide range of constraints, we discuss how to adapt our architecture network to be able to handle both discrete and continuous feasible domains, broadening the applicability of our method to various problem settings.

To show the effectiveness and versatility of our method, we focus on three applications: multidimensional knapsack with continuous variables, combinatorial multidimensional knapsack with discrete variables, and an inventory management problem. The self-supervised learning method allows us to recover over 50\% feasible solutions. By tuning the hyper-parameters in the loss function, our method guarantees a high feasibility level, while maintaining a reasonable sub-optimality threshold. Furthermore, our method achieves significant computational speed-ups compared to traditional solvers, where we observe a minimum speed-up of 10 orders of magnitude when solving a single instance using the learned solver compared to classical solvers, emphasizing the practicality of our method for real-world optimization tasks.

\section{Robust Constrained Optimization: Problem Setting and Self-Supervised Learning Method}\label{sec:method}
\vspace{-0.1cm}
We consider the following family of robust optimization (RO) problems indexed by the instance parameter $z \in \calZ \subseteq \bbbr^{d_z}$:
\begin{equation}
	\label{eq:RO}
    x^*(z) =  \begin{array}[t]{ll} 
    \argmax_{x \in \calX} &
     f_z(x) \\
    \sbt &
     g_z(x,u) \le 0 \quad \forall u \in \calU(z)
    \end{array}
\end{equation}
Here, $\calX \subseteq \bbbr^{d_x}$ refers to the domain of the decision variables $x$, and $\calU \subseteq \bbbr^{d_u}$ refers to the domain of the uncertain parameters $u$.
The objective function $f_z: \calX \to \bbbr$ depends on the instance parameter $z$, and is known with certainty given $z$. The constraint function $g_z: \calX \times \calU \to \bbbr^m$, where $m$ refers to the number of nominal constraints, depends on both the parameter $z$ and on the uncertain parameters $u \in \calU$. Let $g_z^j$ refer to the $j^{\text{th}}$ constraint function for $j \in \{1, \ldots, m\}$. Without loss of generality, we assume that the uncertain parameters $u$ occur in the constraints only. 

Throughout the paper, we consider a more structured special case of the family of problem instances \eqref{eq:RO}, which we refer to as ``nominal-parameterized'' instances.
Specifically, in this special case, each $z$ encodes a vector of uncertain parameters, $\hat{u}(z)$, which we refer to as the nominal values of the uncertain parameters.
The version of \eqref{eq:RO} without robustness, i.e., with only constraints $g_z(x, \hat{u}(z)) \leq 0$ and $x \in \calX$ is thought of as the nominal version of the problem.
We then construct a norm-based uncertainty set based on considering all perturbations of the nominal values $\hat{u}(z)$ within a ball. Possible choices of the norm $\|\cdot\|$ include the $\ell_\infty$ norm corresponding to a box uncertainty set, and a quadratic norm corresponding to an ellipsoidal uncertainty set. This special class of problem instances is formalized in the below assumption, along with an additional assumption ensuring feasibility of \eqref{eq:RO}.

\begin{assumption}\label{assump:1}
The family of RO problems \eqref{eq:RO}, parameterized by $z \in \calZ$, satisfies:
\begin{enumerate}
    \item For some fixed $\rho > 0$, there exists a mapping $\hat{u} : \calZ \to \calU$ to ``nominal uncertain parameters'' and $\calU(z) = \{u : \|u - \hat{u}(z)\| \leq \rho\}$.
    \item For all $z \in \calZ$, the robust feasible set $\{x \in \calX : g_z(x,u) \le 0 \ \forall u \in \calU(z)\}$ is non-empty.
\end{enumerate}
\end{assumption} 
It is important to emphasize that we do not make any strong assumptions on the domain of decision variables $\calX$, and our methodology allows for both convex and integer domains.
These conditions in Assumption \ref{assump:1} are not overly restrictive and are standard in robust optimization. 

\paragraph{{\bf Learning Methodology and Self-Supervised Loss.}}
Our methodology falls under the framework of ``learning to optimize,'' whereby our goal is to approximate the mapping $z \mapsto x^\ast(z)$, specified by the family of RO problems \eqref{eq:RO}, using machine learning. 
We use $h: \calZ \rightarrow \calX$ to denote such an approximation function, referred to as an optimization proxy/surrogate.
The function $h$ is selected from the hypothesis class $\calH$ of candidate functions, which in all of our examples will be a class of neural network models. It is important that $h$ returns outputs that are guaranteed to lie in the domain $\calX$, and we elaborate on how to ensure this later in this paper.
We consider a self-supervised setting whereby we have available a dataset $D= \{z_1, z_2,\ldots,z_n\}$ of realizations of instances of \eqref{eq:RO}. Notably, we do not assume availability of $x^\ast(z_i)$, hence our methodology is {\em self-supervised}. Finally, we also assume availability of a loss function $L : \calX \times \calZ \to \bbbr$, whereby $L(x, z)$ measures the quality of the solution $x$ for the instance of \eqref{eq:RO} parameterized by $z$. Putting all of these ingredients together and applying the empirical risk minimization (ERM) principle of machine learning (a.k.a. sample average approximation) leads to the following training problem:
\begin{equation}\label{eq:ERM}
\begin{aligned}
    \min_{h \in \calH } \ \frac{1}{n}\sum_{i=1}^n L (h(z_i),z_i)
\end{aligned}
\end{equation}
The default supervised learning approach would treat this setting as a regression problem and choose a regression based loss function such as $L^{\mathrm{SL}}(h(z),z) = \|h(z) - x^\ast(z)\|_2^2$. Instead of adopting a supervised learning approach, which requires access to pre-solved instances $x^\ast(z_i)$, we employ a penalty-based self-supervised loss function that directly encodes the problem structure by balancing optimality and feasibility with respect to the robust constraints. This approach is particularly useful in our setting where obtaining pre-solved instances can be computationally expensive due to the increased complexity of solving robust optimization counterparts over their nominal versions \cite{ben2009robust,ben2015deriving}.
Inspired by the exact penalty method used to solve constrained optimization \cite{di1989exact}, we define our self-supervised loss function $L_{\nu}^{\mathrm{SSL}} : \calX \times \calZ \to \bbbr$, given a penalty parameter $\nu > 0$, as follows:
\begin{equation}\label{eq:ssl_loss}
L_{\nu}^{\mathrm{SSL}}(x, z) := -f_z(x) + \nu\sum_{j = 1}^m \left[\max_{u \in \calU(z)} g_z^j(x,u)\right]^+,
\end{equation}
where $[\cdot]^+ := \max\{0, \cdot\}$ is the exact penalty function. 
The first term of the self-supervised loss function aims to maximize the objective function $f(\cdot)$, while the second term is a penalty that attributes a positive weight to any unsatisfied robust constraints.
\begin{assumption}\label{assump:2}
The learning problem \eqref{eq:ERM} is practically tractable in the sense:
\begin{enumerate} 
    \item For all $ z \in \calZ $, the objective function $ f_z(\cdot) $ is continuous and differentiable almost everywhere with respect to $x$ on an open set containing $\calX$ .
    \item For each $j \in \{1, \ldots, m\}$, the function $\bar g_z^j(x) =  \left[\max_{u \in \calU(z)} g_z^j(x,u)\right]^+$ is continuous and differentiable almost everywhere with respect to $x$ on an open set containing $\calX$ .
\end{enumerate} 
\end{assumption}
Assumption \ref{assump:2} ensures practical tractability of our learning problem when using a neural network hypothesis class, since the required differentiability properties ensure that we can apply modern automatic differentiation frameworks such as PyTorch \cite{paszke2019pytorch}.
While part (1.) of the Assumption \ref{assump:2} is readily satisfied, part (2.), however, requires that the robust counterpart of the RO problem \eqref{eq:RO} has a tractable reformulation with differentiable constraints. 
Under the assumption that, for all $j \in \{1, \ldots, m\}$ and all fixed $x \in \calX$, the constraint function $g_z^j(x,\cdot)$ is concave in the uncertain parameters $u$, the results of Ben-Tal \cite{ben2015deriving} provide a framework for constructing such a tractable robust counterpart for a wide range of uncertainty sets and constraint functions.

\paragraph{{\bf Ensuring Domain Feasibility via Neural Network Structure.}}
While the previously introduced self-supervised loss function is used to promote feasibility w.r.t. the robust inequality constraints, we adopt a different approach to ensure the variable domain constraint $ x \in \calX $. Indeed, the domain constraint $ x \in \calX $ is usually relatively much simpler, including simple continuous or discrete sets, than the robust inequality constraints. As such, specifically in the case of neural network hypothesis classes $\calH$, we discuss how to directly engineer the network to automatically ensure that $h(z) \in \calX$, for several important domains $\calX$.
For many typical convex domains, using a standard activation function at the last layer is sufficient to ensure feasibility of the output. For example, the ReLU activation can be used to ensure that $h(z) \in \calX = \bbbr_{+}^n$ and the sigmoid activation can be used to ensure that $h(z) \in \calX = [0,1]^n$.

The case of an mixed integer or discrete domain set $\calX$ is more challenging. Related to this case is the problem of differentiating through a mixed integer program (MIP), which occurs in decision-focused learning and related areas \cite{wilder2018meldingdatadecisionspipelinedecisionfocused,ferber2019mipaalmixedintegerprogram,vlastelica2020differentiationblackboxcombinatorialsolvers}. Most of the literature addressing this or related challenges develop gradient-based approaches that are suited for differentiating through a generic MIP with complex constraints. In our case, we implicitly assume that the domain constraints represented by $\calX$ are relatively simple and that the ``complex'' constraints are the robust ones that have been incorporated into the self-supervised loss function.
As such, we adopt an approach based on the assumption that we can efficiently solve linear optimization problems over $\calX$ and also a related family of ``smoothed'' quadratic optimization problems that enhance differentiability properties \cite{wilder2018meldingdatadecisionspipelinedecisionfocused}. Our approach is inspired by \cite{ferber2023surcolearninglinearsurrogates}, who propose using a linear integer optimization surrogate to address learning the solutions of a family of mixed integer nonlinear optimization problems, as well as the older literature on structured prediction methods \cite{bakir2007predicting,osokin2017structured}

To be precise, let us assume that $\calX \subseteq \bbbr^{d_x}$ is a non-empty compact mixed integer set. To ensure feasibility $h(z) \in \calX$, we carefully design the last layer of our neural network models $h \in \calH$. The previous layers can consist of any standard neural network architecture. The last layer should take in as input a vector $w \in \bbbr^{d_x}$ and output one of two deterministic mappings:  (i) $g_{\calX}^{\mathrm{train}}(w)$ at training time, and (ii) $g_{\calX}^{\mathrm{test}}(w)$ at testing/inference time. Let us first describe $g_{\calX}^{\mathrm{test}}$. We let $g_{\calX}^{\mathrm{test}}(w) \in \argmax_{x \in \calX} w^\top x$ be an arbitrary solution to a linear optimization problem over $\calX$, given the cost vector $w \in \bbbr^{d_x}$ determined by the computations at the earlier layers. Thus, one of our implicit assumptions about the simplicity of $\calX$ is that it is tractable to solve linear optimization problems. We design $g_{\calX}^{\mathrm{train}}$ as an approximation of $g_{\calX}^{\mathrm{test}}$, with better differentiability properties. Specifically, let $\bar{\calX} = \mathrm{conv}(\calX)$ denote the convex hull of $\calX$. Given some smoothing parameter $\gamma > 0$, we propose the following smooth 
$g_{\calX}^{\mathrm{train}}(w) = \argmax_{x \in \bar{\calX}} \left\{ w^\top x - \frac{\gamma}{2} \|x\|_2^2 \right\}. $

The second implicit assumption regarding $\calX$ is that the above optimization problem is also readily solvable. Note that the above formulation is a continuous optimization problem over the convex hull $\bar{\calX}$, and the mapping can be differentiated w.r.t. $w$ by implicit differentiation or other techniques \cite{wilder2018meldingdatadecisionspipelinedecisionfocused}. Inspired by \cite{qi2023iceo}, we further propose a learning-based approximation to $g_{\calX}^{\mathrm{train}}(w)$ to enhance the smoothness properties and ease of computation of evaluating the gradient. This approximation leads to a 2-phase approach, where the first phase trains a simple network approximating $g_{\calX}^{\mathrm{train}}(w)$ using a regression-based model. Once this model is trained, we freeze the weight of this network and use it as the last layer in the optimization predictive model, which we then train using our self-supervised approach.
\begin{example}
As a concrete example, consider the case where $x \in \calX = \{0, 1\}^{d_x}$ corresponds to binary integer constraints. To ensure feasibility with respect to $\calX$ at testing time, we note that $g_{\calX}^{\mathrm{test}}(w) \in \argmax_{x \in \calX} w^\top x$ is given by the component-wise step function $[g_{\calX}^{\mathrm{test}}(w)]_i = \mathbf{1}(w_i > 0)$ where $\mathbf{1}(\cdot)$ is an indicator function equal to 1 if the argument is true and 0 otherwise.

Similarly the approximation function $g_{\calX}^{\mathrm{train}}(w)$
can also be expressed component-wise as $[g_{\calX}^{\mathrm{train}}(w)]_i = \min\left(1,\max\left(0,w_i/\gamma\right)\right)$. The analytic expression can either be directly differentiated or approximated using a neural network, ensuring smoothness and facilitating efficient optimization through gradient-based methods.
\end{example}
\section{Experimental Results}
\vspace{-0.1cm}
In this section, we present the numerical evaluation of our method, focusing on the following applications:
\begin{application}[Multidimensional Knapsack Problem]
We consider a multidimensional knapsack problem where the goal is to allocate a set of items, with a vector of values $ v_z $ and a matrix of weights $ W $, across multiple knapsacks with capacities $ C_z $. The items' weight matrix $ W $ belongs to a known uncertainty set $ \mathcal{U}(z) $. The problem can be formulated as follows:
\begin{equation*}
    \begin{array}{llll}
    x^*(z) =& \max\limits_{x \in \mathcal{X}} & v_z^\top x & \\
    & \text{s.t.} & W x \leq C_z, & \quad \forall W \in \mathcal{U}(z).
    \end{array}
\end{equation*}
We consider both the discrete case where $ \mathcal{X} = \{0,1\}^{d_x} $ and the continuous case where $ \mathcal{X} = [0,1]^{d_x} $. We generate a synthetic dataset $\mathcal{D} = \{(v_i, \hat{W}_i, C_i)\}_{i=1}^{N}$ for the knapsack problem, where $ \hat{W}_i $ represents the nominal value of the weight matrix.
\end{application}

\begin{application}[Inventory Management Problem] We consider an inventory management problem inspired by the classical newsvendor model, extended to a multi-retailer system with a centralized warehouse, with a formulation adopted from \cite{iancu2014pareto}. $N$ retailers face uncertain demand $d_z = d_z^0 + Q_z u$, where $d_z^0$ is the expected demand,  
$u$ is a $k$-dimensional vector of uncertainty factors, and 
$Q_z$ is a matrix capturing the sensitivity of retailers' demand to these factors. The stocking decisions $x$ aim to maximize profit across all retailers. Considering an auxiliary second stage variable $y(u)$ representing sales (as in \cite{iancu2014pareto}), the problem is formulated as the following two-stage adjustable robust optimization:
\begin{equation*}
    \begin{array}{llc}
    \max_{x \in \calX, P \in \bbbr} &P \\
    \text{s.t.} & P \le  r_z^\top y(u)  - {c_z^o}^\top x & \forall u \in \calU(z) \\
    & y(u) = \min(x,d_z^0 + Q_z u)  & \forall u \in \calU(z) \\
            & \mathbf{1}^\top x \leq C_z
    \end{array}
\end{equation*}
where the vector $x$ is contained in $\calX = \{x : 0\le x \le c\}$, and $r_z$ and $c_z^o$ are the revenue and cost  per unit sold respectively. Using a linear decision rule \cite{bertsimas2022robustadaptiveoptim} for the second-stage decisions, i.e., $y(u)= Y  u + y_0$, we derive the following one-stage robust problem:
\begin{equation*}
\begin{array}{llc}
\max_{P \in \bbbr, Y  \in \bbbr^{N \times k}, y_0 \in \bbbr^N, x \in \mathcal{X}} & P &\\
\text{s.t.} & P \leq r_z^\top (Y u +y_0) - {c_z^0}^\top x & \forall u \in \mathcal{U}_z \\
            & Y u +y_0 \leq x, & \forall u \in \mathcal{U}_z \\
            & Y u +y_0 \leq d_z^0 + Q_z u, & \forall u \in \mathcal{U}_z \\
            & \mathbf{1}^\top x \leq C_z.
\end{array}
\end{equation*}
We generate a synthetic dataset $\mathcal{D} = \{(r_i, c_i^0,d^0_i,Q_i,\hat{u}_i)\}_{i=1}^{N}$ for the inventory management problem, where $\hat{u}_i $ represents the nominal value of the uncertainty factor.
\end{application}

For each application, we use a synthetically generated dataset to train a two-layer fully connected neural network with problem-specific activation functions, as described in Section \ref{sec:method}, to approximate the optimal solution of the parametric optimization problem. The dataset is split into training (70\%), validation (15\%), and testing (15\%). At test time, we compute the optimal robust solution using the Gurobi solver and evaluate the performance of our method based on feasibility, optimality, and computational efficiency. To assess feasibility, we compute the average maximum constraint violation over the test set, defined as $\max_{j} \left[ \bar{g}_z^j\left(h(z)\right) \right]$, where $\bar{g}_z^j$ is given in Assumption \ref{assump:2}, and report the percentage of feasible solutions produced by the learned model. For optimality, we measure regret defined as $(f_z^*-\hat{f}_z)/f_z^*$, where $f_z^*$ is the optimal objective value obtained by Gurobi, and $\hat f_z$ is the objective value of the solution provided by the learned model. Finally, to evaluate computational efficiency, we compare the average computational time required by both Gurobi and the learned solver reported in seconds.

\begin{table}[h]
\caption{Experimental Results for Constraint Optimization Problems}
\label{tab:experimental_results}
\resizebox{\textwidth}{!}{
\begin{tabular}{l|c|c|c|c|c|c|c}
\toprule
\textbf{Application} & \textbf{Problem} & \textbf{Model \&}         & \textbf{Optimality} & \textbf{Max Constraint}  & \textbf{\% Feasible} & \textbf{Inference}  & \textbf{Solver} \\ 
                     & \textbf{Size}    & \textbf{Penalty Coefficient} & \textbf{Regret}       & \textbf{Violation}  & \textbf{Instances}   & \textbf{time} & \textbf{time}\\ 
\midrule
Knapsack & $d_x= 50$  & SL / & -1.41 & 3.93 & 0.88 & \textbf{0.0023} & 0.0279 
\\ & $m=5$ & SSL /1 & \textbf{0.4} & 0.36& 57.6 &  &
\\ &  & SSL /10 & 0.64 & 0.082& 95.6 &  &
\\ &  & SSL /20 & 0.63 & 0.088 & 98 &  &
\\ &  & SSL /50 &0.70&\textbf{ 0.033} & \textbf{99} &  &\\
\midrule
IP Knapsack & $d_x= 50$ & SL / & -0.74 & 3.82 & 2 & \textbf{0.0027 }& 0.1311
\\ & $m=5$  & SSL /1 & 0.50 & 0.44 &  60 &  &
\\ &  & SSL /10 & 0.79 & \textbf{0.06} & \textbf{94} &  &
\\ &  & SSL /20 & \textbf{0.22 }& 0.33 & 84 &  &
\\ &  & SSL /50 & 0.30 & 0.38 & 59 &  &
\\
\midrule
Inventory & $d_x=N=50$   & SL / & \textbf{-0.008} & 90.9 & 0 &  \textbf{0.0011 }& 0.3975
\\ & $d_u=k=5$ & SSL /50& -1.8 & 0.065 & 22.6&  &
\\ &  & SSL /100 & -2.13 & 0.112 & 57.5&  &
\\ &  & SSL /200 &  -1.94 & 0.0150 & 72.4&  &
\\ &  & SSL /500 & -2.21 &\textbf{ 0.012} &  \textbf{91.7}&  &
\\
\bottomrule
\end{tabular}
}
\end{table}

Compared to the supervised learning approach, the solution provided by our method achieves a high level of feasibility, as shown in Table \ref{tab:experimental_results}, whereby the percentage of feasible solutions returned by the self-supervised learning approach is consistently high across different applications and penalty coefficients. This highlights the effectiveness of our approach in generating feasible solutions, addressing a key limitation of supervised learning methods, which often produce infeasible solutions.  
Table \ref{tab:experimental_results} also highlights the trade-off between feasibility and optimality depending on the choice of the penalty coefficient used during training. Exploring this trade-off may be valuable in real-world applications, where the relative importance of feasibility and optimality depends on the specific requirements of the problem. Additionally, Table \ref{tab:experimental_results} demonstrates the computational efficiency of our method, achieving a speedup of over 10 orders of magnitude compared to a standard optimization solver like Gurobi. This efficiency can further be enhanced by the learned solver's ability to perform batch predictions, which may be extremely valuable when solving multiple problems simultaneously.
\section{Conclusion}
\vspace{-0.1cm}
We present a self-supervised learning approach to solve optimization problems under uncertainty leveraging tractable reformulations of robust optimization problems. Our approach relies on penalty based self-supervised learning loss function to learn robust solutions, eliminating the need for pre-solved instances. Additionally, this work utilizes neural network models capable of outputting values for both integer and discrete variables, allowing our framework to handle a broad range of optimization problems. Numerical experiments validate the effectiveness of the proposed approach, revealing a trade-off between the feasibility and optimality of the solutions. Additionally, our method achieves significant computational speedup compared to traditional solvers which can be crucial in real-time decision-making applications.
\begin{credits}

\subsubsection{\ackname} 
This research was supported by NSF AI Institute for Advances in Optimization Award 2112533.

\end{credits}

\bibliographystyle{splncs04}
\bibliography{main}

\end{document}